# A NOTE ON THE GUROV-RESHETNYAK CONDITION

A.A. KORENOVSKYY, A.K. LERNER, AND A.M. STOKOLOS

ABSTRACT. An equivalence between the Gurov-Reshetnyak $GR(\varepsilon)$ and Muckenhoupt $A_\infty$ conditions is established. Our proof is extremely simple and works for arbitrary absolutely continuous measures.

Throughout the paper, $\mu$ will be a positive measure on $\mathbb{R}^n$ absolutely continuous with respect to Lebesgue measure. Denote

$$\Omega_\mu(f;Q) = \frac{1}{\mu(Q)} \int_Q |f(x) - f_{Q,\mu}| d\mu(x), \quad f_{Q,\mu} = \frac{1}{\mu(Q)} \int_Q f(x) d\mu(x).$$

**Definition 1.** We say that a nonnegative function $f$, $\mu$-integrable on a cube $Q_0$, satisfies the Gurov-Reshetnyak condition $GR_\mu(\varepsilon)$, $0 < \varepsilon < 2$, if for any cube $Q \subset Q_0$,

(1) $$\Omega_\mu(f;Q) \leq \varepsilon f_{Q,\mu}.$$

When $\mu$ is Lebesgue measure we drop the subscript $\mu$.

This condition appeared in [6, 7]. It is important in Quasi-Conformal Mappings, PDEs, Reverse Hölder Inequality Theory, etc. (see, e.g., [2, 8]). Since (1) trivially holds for all positive $f \in L_\mu(Q_0)$ if $\varepsilon = 2$, only the case $0 < \varepsilon < 2$ is of interest. It was established in [2, 6, 7, 8, 13] for Lebesgue measure and in [4, 5] for doubling measures that if $\varepsilon$ is small enough, namely $0 < \varepsilon < c2^{-n}$, the $GR_\mu(\varepsilon)$ implies $f \in L^p_\mu(Q_0)$ with some $p > 1$ depending on $\varepsilon$. The machinery used in the articles mentioned above does not work for $\varepsilon > 1/8$ even in the one-dimensional case.

The one-dimensional improvement of these results was done in [9]. Namely, for any $0 < \varepsilon < 2$ it was proved that $GR(\varepsilon) \subset L^p_{loc}$ where $1 < p < p(\varepsilon)$; moreover a sharp bound $p(\varepsilon)$ for the exponent was discovered. The main tool in [9] is the Riesz Sunrising Lemma which has no multidimensional version since it involves the structure of open sets on a real line.

In the present article using simple arguments we prove that for any $n \geq 1$, $0 < \varepsilon < 2$ and arbitrary absolutely continuous measure $\mu$, the Gurov-Reshetnyak condition $GR_\mu(\varepsilon)$ implies the weighted $A_\infty(\mu)$ Muckenhoupt condition. And conversely, $A_\infty(\mu)$ implies $GR_\mu(\varepsilon_0)$ for some $0 < \varepsilon_0 < 2$.

In the non-weighted (or doubling) case R.R. Coifman and C. Fefferman [3] have found several equivalent descriptions of the $A_\infty$ property. Recently these descriptions have been transferred by J. Orobitg and C. Perez [12] to the non-doubling case. For our purposes it will be convenient to define $A_\infty(\mu)$ by the following way.

**Definition 2.** We say that a nonnegative function $f$, $\mu$-integrable on a cube $Q_0$, satisfies Muckenhoupt condition $A_\infty(\mu)$ if there exist $0 < \alpha, \beta < 1$ such that for any cube $Q \subset Q_0$,

$$\mu\{x \in Q : f(x) > \beta f_{Q,\mu}\} > \alpha \mu(Q).$$

Our main result is the following

1991 *Mathematics Subject Classification*. Primary 42B25.

The first author was partially supported by Ukrainian Foundation of Fundamental Research, grant F7/329 - 2001.

**Theorem 1.**

(i) *Suppose that for some $0 < \varepsilon < 2$ nonnegative function $f$ satisfies the inequality*
$$\Omega_\mu(f;Q) \leq \varepsilon f_{Q,\mu}.$$

*Then, for $\varepsilon < \lambda < 2$ we have*
$$\mu\{x \in Q : f(x) > (1 - \varepsilon/\lambda) f_{Q,\mu}\} \geq (1 - \lambda/2)\mu(Q);$$

(ii) *Suppose that for some $0 < \alpha, \beta < 1$ nonnegative function $f$ satisfies the inequality*
$$\mu\{x \in Q : f(x) > \beta f_{Q,\mu}\} > \alpha\mu(Q).$$

*Then*
$$\Omega_\mu(f;Q) \leq 2(1 - \alpha\beta) f_{Q,\mu}.$$

*Proof.* Set $E = \{x \in Q : f(x) > (1 - \varepsilon/\lambda) f_{Q,\mu}\}$, $E^c = Q \setminus E$. Suppose that $\mu(E^c) > 0$ (otherwise part (i) is trivial). Then

$$\frac{\varepsilon}{\lambda} f_{Q,\mu} \leq \inf_{x \in E^c}\bigl(f_{Q,\mu} - f(x)\bigr) \leq \frac{1}{\mu(E^c)} \int_{E^c} \bigl(f_{Q,\mu} - f(x)\bigr) d\mu(x)$$

$$\leq \frac{1}{\mu(E^c)} \int_{\{x \in Q : f(x) < f_{Q,\mu}\}} \bigl(f_{Q,\mu} - f(x)\bigr) d\mu(x) = \frac{1}{\mu(E^c)} \frac{\mu(Q)}{2} \Omega_\mu(f;Q) \leq \frac{1}{\mu(E^c)} \frac{\mu(Q)}{2} \varepsilon f_{Q,\mu}.$$

Hence, $\mu(E^c) \leq (\lambda/2)\mu(Q)$, as required.

To prove the second part, set $E = \{x \in Q : f(x) > \beta f_{Q,\mu}\}$, $E^c = Q \setminus E$. Then $\mu(E^c) \leq (1 - \alpha)\mu(Q)$, and we have

$$\Omega_\mu(f;Q) = \frac{2}{\mu(Q)} \int_{\{x \in Q: f(x) < f_{Q,\mu}\}} \bigl(f_{Q,\mu} - f(x)\bigr) d\mu(x)$$

$$= \frac{2}{\mu(Q)} \int_{\{x \in Q : \beta f_{Q,\mu} < f(x) < f_{Q,\mu}\}} \bigl(f_{Q,\mu} - f(x)\bigr) d\mu(x) + \frac{2}{\mu(Q)} \int_{E^c} \bigl(f_{Q,\mu} - f(x)\bigr) d\mu(x)$$

$$\leq \frac{2}{\mu(Q)} f_{Q,\mu} \Bigl((1 - \beta)\mu(E) + \mu(E^c)\Bigr) = \frac{2}{\mu(Q)} f_{Q,\mu} \Bigl((1 - \beta)\mu(Q) + \beta\mu(E^c)\Bigr)$$

$$\leq \frac{2}{\mu(Q)} f_{Q,\mu} \Bigl((1 - \beta)\mu(Q) + \beta(1 - \alpha)\mu(Q)\Bigr) = 2(1 - \alpha\beta) f_{Q,\mu}.$$

The theorem is proved. □

**Corollary.** *The following characterization of $A_\infty(\mu)$ holds:*
$$A_\infty(\mu) = \bigcup_{0 < \varepsilon < 2} GR_\mu(\varepsilon).$$

Since $A_\infty(\mu)$ condition is equivalent to the weighted reverse Hölder inequality for some $p > 1$ (cf. [12]), i.e.

(2) $$\Bigl(\frac{1}{\mu(Q)} \int_Q (f(x))^p d\mu(x)\Bigr)^{1/p} \leq c \frac{1}{\mu(Q)} \int_Q f(x) d\mu(x) \quad (Q \subset Q_0),$$

A NOTE ON THE GUROV-RESHETNYAK CONDITION 3we see that for any $0 < \varepsilon < 2$ a function $f$ satisfying the Gurov-Reshetnyak condition $GR_\mu(\varepsilon)$ belongs to $L^p_\mu(Q_0)$ for some $p > 1$. Observe that such approach (i.e. $GR_\mu(\varepsilon) \Rightarrow A_\infty(\mu) \Rightarrow$ Reverse Hölder) does not give the optimal order of integrability for small $\varepsilon$, though it is known [2, 4] in doubling case that $GR_\mu(\varepsilon) \subset L^{p(\varepsilon)}_\mu(Q_0)$, where $p(\varepsilon) \asymp c_n/\varepsilon, \varepsilon \to 0$, and this order is sharp. However we will show that part (i) of Theorem 1 allows us to obtain the same order for any measure $\mu$, any $0 < \varepsilon < 2$, and $f \in GR_\mu(\varepsilon)$. We will need the following

**Covering Lemma** [10]. *Let $E$ be a subset of $Q_0$, and suppose that $\mu(E) \leq \rho\mu(Q_0)$, $0 < \rho < 1$. Then there exists a sequence $\{Q_i\}$ of cubes contained in $Q_0$ such that*

(i) $\mu(Q_i \cap E) = \rho\mu(Q_i)$;
(ii) *the family $\{Q_i\}$ is almost disjoint with constant $B(n)$, that is, every point of $Q_0$ belongs to at most $B(n)$ cubes $Q_i$;*
(iii) $E' \subset \cup_j Q_j$, *where $E'$ is the set of $\mu$-density points of $E$.*

Recall that the non-increasing rearrangement of $f$ on a cube $Q_0$ with respect to $\mu$ is defined by

$$f^*_\mu(t) = \sup_{E \subset Q_0 : \mu(E) = t} \inf_{x \in E} |f(x)| \quad (0 < t < \mu(Q_0)).$$

Denote $f^{**}_\mu(t) = t^{-1} \int_0^t f^*_\mu(\tau) d\tau$.

**Theorem 2.** *Let $0 < \varepsilon < 2$, and $f \in GR_\mu(\varepsilon)$. Then for $\varepsilon < \lambda < 2$, $\rho < 1 - \lambda/2$, and $t \leq \rho\mu(Q_0)$ we have*

$$(3) \qquad f^{**}_\mu(t) \leq \left(B(n)\frac{\lambda/\rho + 1}{\lambda - \varepsilon}\varepsilon + 1\right)f^*_\mu(t).$$

**Remark.** A well-known argument due to Muckenhoupt (see, e.g., [11, Lemma 4]) shows that (3) implies the reverse Hölder inequality (2) for all $p < 1 + \left(\frac{\lambda - \varepsilon}{B(n)(\lambda/\rho + 1)}\right)\frac{1}{\varepsilon}$.

*Proof of Theorem 2.* Set $E = \{x \in Q_0 : f(x) > f^*_\mu(t)\}$, and apply the Covering Lemma to $E$ and number $\rho$. We get cubes $Q_i \subset Q_0$, satisfying (i)-(iii). Since $\rho < 1 - \lambda/2$, we obtain from (i) that for each $Q_i$,

$$(f\chi_{Q_i})^*_\mu((1 - \lambda/2)\mu(Q_i)) \leq f^*_\mu(t).$$

Hence, by Theorem 1,

$$(4) \qquad f_{Q_i,\mu} \leq \frac{\lambda}{\lambda - \varepsilon}(f\chi_{Q_i})^*_\mu((1 - \lambda/2)\mu(Q_i)) \leq \frac{\lambda}{\lambda - \varepsilon}f^*_\mu(t),$$

and so,

$$(5) \qquad \Omega_\mu(f; Q_i) \leq \frac{\varepsilon\lambda}{\lambda - \varepsilon}f^*_\mu(t).$$

Further, by (ii),

$$\sum_i \mu(Q_i \cap E) \leq B(n)\mu(E) \leq B(n)t.$$

Therefore, using a well-known property of rearrangement (see, e.g. [1]) and (4), (5), we obtain

$$\begin{aligned}
t\big(f_\mu^{**}(t) - f_\mu^*(t)\big) &= \int_E \big(f(x) - f_\mu^*(t)\big)d\mu(x) = \sum_i \int_{E \cap Q_i} \big(f(x) - f_\mu^*(t)\big)d\mu(x) \\
&= \sum_i \int_{E \cap Q_i} \big(f(x) - f_{Q_i,\mu}\big)d\mu(x) + \sum_i \mu(E \cap Q_i)\big(f_{Q_i,\mu} - f_\mu^*(t)\big) \\
&\leq \frac{\varepsilon \lambda}{\lambda - \varepsilon} f_\mu^*(t) \sum_i \mu(Q_i) + \frac{\varepsilon}{\lambda - \varepsilon} f_\mu^*(t) \sum_i \mu(E \cap Q_i) \\
&\leq B(n) \frac{\lambda/\rho + 1}{\lambda - \varepsilon} \varepsilon t f_\mu^*(t),
\end{aligned}$$

which gives the desired result. $\square$

## References

[1] C. Bennett and R. Sharpley, *Interpolation of operators*, Academic Press, New York, 1988.
[2] B. Bojarski, *Remarks on the stability of reverse Hölder inequalities and quasi-conformal mappings*, Ann. Acad. Sci. Fenn. Ser. A I Math. **10** (1985), 89–94.
[3] R.R. Coifman and C. Fefferman, *Weighted norm inequalities for maximal functions and singular integrals*, Studia Math. **15** (1974), 241–250.
[4] M. Franciosi, *Weighted rearrangements and higher integrability results*, Studia Math. **92** (1989), no. 2, 131–139.
[5] M. Franciosi, *On weak reverse integral inequalities for mean oscillations*, Proc. Amer. Math. Soc. **113** (1991), no. 1, 105–112.
[6] L.G. Gurov, *On the stability of Lorentz mappings. Estimates for derivatives* (Russian, English), Sov. Math., Dokl. **16** (1975), 56-59; translation from Dokl. Akad. Nauk SSSR **220** (1975), 273-276.
[7] L.G. Gurov and Yu.G. Reshetnyak, *A certain analogue of the concept of a function with bounded mean oscillation* (Russian), Sibirsk. Mat. Zh. **17** (1976), no. 3, 540–546.
[8] T. Iwaniec, *On $L^p$-integrability in PDEs and quasiregular mappings for large exponents*, Ann. Acad. Sci. Fenn. Ser. A I Math. **7** (1982), no. 2, 301–322.
[9] A.A. Korenovskii, *The connection between mean oscillations and exact exponents of summability of functions* (Russian), Mat. Sb. **181** (1990), no. 12, 1721–1727; translation in Math. USSR-Sb. **71** (1992), no. 2, 561–567
[10] J. Mateu, P. Mattila, A. Nicolau and J. Orobitg, *BMO for nondoubling measures*, Duke Math. J. **102**, No.3, (2000), 533-565.
[11] B. Muckenhoupt, *Weighted inequalities for the Hardy maximal function*, Trans. Amer. Math. Soc. **165** (1972), 207-227.
[12] J. Orobitg, and C. Perez, *$A_p$ weights for nondoubling measures in $R^n$ and applications*, Trans. Amer. Math. Soc. **354** (2002), no. 5, 2013–2033.
[13] I. Wik, *Note on a theorem by Reshetnyak-Gurov*, Studia Math. **86** (1987), no. 3, 287–290.

Department of Mathematical Analysis, IMEM, National University of Odessa, Dvoryanskaya, 2, 65026 Odessa, Ukraine
*E-mail address*: anakor@paco.net

Department of Mathematics and Computer Science, Bar-Ilan University, 52900 Ramat Gan, Israel
*E-mail address*: aklerner@netvision.net.il

Department of Mathematics, University of Connecticut, U-9, Storrs, CT 06268, USA
*E-mail address*: stokolos@math.uconn.edu